\documentclass{amsart}
\usepackage{amsmath,amsfonts,amsthm,amssymb,epsfig,graphicx,eucal,natbib,dsfont,accents,url,epstopdf}

\DeclareFontFamily{U}{mathx}{\hyphenchar\font45}
\DeclareFontShape{U}{mathx}{m}{n}{
      <5> <6> <7> <8> <9> <10>
      <10.95> <12> <14.4> <17.28> <20.74> <24.88>
      mathx10
      }{}
\DeclareSymbolFont{mathx}{U}{mathx}{m}{n}
\DeclareFontSubstitution{U}{mathx}{m}{n}
\DeclareMathAccent{\widecheck}{0}{mathx}{"71}

\newtheorem{thm}{Theorem}[section]

\numberwithin{figure}{section}
\numberwithin{table}{section}

\theoremstyle{definition}

\theoremstyle{remark}
\newtheorem{rem}[thm]{Remark}

\newcommand*{\dt}[1]{%
  \accentset{\mbox{\large\bfseries .}}{#1}}
\newcommand{\dd}{\mathrm{d}}
\DeclareMathAlphabet{\mathbfit}{OML}{cmm}{b}{it}

\begin{document}
\title[GOF Testing for Copulas: A Distribution-Free Approach]{Goodness-of-Fit Testing for Copulas: A Distribution-Free Approach}
\author[S. U. Can]{Sami Umut Can}
\address{Department of Quantitative Economics, University of Amsterdam, P.O. Box 15867, 1001 NJ Amsterdam, The Netherlands}
\email{s.u.can@uva.nl}
\author[J. H. J. Einmahl]{John H.J. Einmahl}
\address{Department of Econometrics \& OR and CentER, Tilburg University, P.O. Box 90153, 5000 LE Tilburg, The Netherlands}
\email{j.h.j.einmahl@tilburguniversity.edu}
\author[R. J. A. Laeven]{Roger J.A. Laeven}
\address{Department of Quantitative Economics, University of Amsterdam, P.O. Box 15867, 1001 NJ Amsterdam, The Netherlands}
\email{r.j.a.laeven@uva.nl}

\begin{abstract}
Consider a random sample from a continuous multivariate distribution function $F$ with copula $C$. In order to test the null hypothesis that $C$ belongs to a certain parametric family, we construct an empirical process on the unit hypercube that converges weakly to a standard Wiener process under the null hypothesis. This process can therefore serve as a `tests generator' for asymptotically distribution-free goodness-of-fit testing of copula families. We also prove maximal sensitivity of this process to contiguous alternatives. Finally, we demonstrate through a Monte Carlo simulation study that our approach has excellent finite-sample performance, and we illustrate its applicability with a data analysis.
\end{abstract}


\maketitle

\section{Introduction}\label{sec.intro}
Consider a $d$-variate ($d \ge 2$) distribution function (df) $F$ with continuous margins $F_1, \ldots, F_d$. By the representation theorem of \cite{sklar:1959}, there is a unique df $C$ on the unit hypercube $[0,1]^d$ with uniform margins such that
\begin{equation}\label{cop}
F(\mathbf{x}) = C(F_1(x_1), \ldots, F_d(x_d)), \quad \mathbf{x} = (x_1, \ldots, x_d)^\mathsf{T} \in \mathbb{R}^d.
\end{equation}
In fact, if $\mathbf{X} = (X_1, \ldots, X_d)^\mathsf{T}$ is a random vector with joint df $F$, then it is easily seen that the unique df $C$ satisfying (\ref{cop}) is the joint df of the com\-ponent-wise probability integral transforms $(F_1(X_1), \ldots, F_d(X_d))^\mathsf{T}$, given by
\begin{equation}\label{cop2}
C(\mathbf{u}) = F(Q_1(u_1), \ldots, Q_d(u_d)), \quad \mathbf{u} = (u_1, \ldots, u_d)^\mathsf{T} \in [0,1]^d,
\end{equation}
with $Q_j$ denoting the left-continuous quantile function associated with $F_j$, i.e.\ $Q_j(\cdot) = \inf\{x \in \mathbb{R} : F_j(x) \ge \cdot\}$, for $j=1, \ldots, d$.

The df $C$ satisfying (\ref{cop}) or (\ref{cop2}) is called the \emph{copula} associated with $F$, and it is a representation of the dependence structure between the margins of $F$, since $C$ contains no information about the margins, yet together with the margins it characterizes $F$. Thus copulas allow separate modeling of margins and dependence structure in multivariate settings, which has proved to be a useful approach in a wide range of applied fields, from medicine and climate research to finance and insurance. We refer to recent comprehensive monographs such as \cite{nelsen:2006}, \cite{joe:2015} and \cite{durante:sempi:2016} for more background on copula theory and its various applications.

The present work is concerned with goodness-of-fit (GOF) testing for copulas. More specifically, we assume that an i.i.d.\ sample
\begin{equation*}
\mathbf{X}_1 = (X_{11}, \ldots, X_{1d})^\mathsf{T}, \ldots, \mathbf{X}_n = (X_{n1}, \ldots, X_{nd})^\mathsf{T}
\end{equation*}
is observed from an unknown $d$-variate df $F$ with continuous margins $F_1, \ldots,$ $F_d$ and copula $C$ as above. We are interested in testing the hypothesis $C \in \mathcal{C}$ against the alternative $C \notin \mathcal{C}$, where $\mathcal{C} = \{C_{\boldsymbol{\lambda}} : \boldsymbol{\lambda} \in \Lambda\}$ denotes a parametric family of copulas, indexed by a finite-dimensional parameter $\boldsymbol{\lambda}$. There is a rich cornucopia of parametric copula families used in various applications (see, e.g., Ch.\ 4 of \cite{joe:2015} or Ch.\ 6 of \cite{durante:sempi:2016} for extensive lists), and new ones are introduced regularly in the literature, so the testing problem just described is clearly very relevant for practitioners making use of copula modeling in their work.

GOF testing for copulas is not a new problem, and several approaches have been proposed in the literature since the early 1990s, each with their advantages and limitations in specific situations. A partial list includes \cite{genest:rivest:1993}, \cite{shih:1998}, \cite{wang:wells:2000}, \cite{breymann:etal:2003}, \cite{fermanian:2005}, \cite{genest:etal:2006} and \cite{dobric:schmid:2007}. There seems to be no single approach that is universally preferred over others. For a broad overview and comparison of various GOF testing procedures for copulas, we refer to \cite{berg:2009}, \cite{genest:etal:2009} and \cite{fermanian:2013}. As pointed out in the latter papers, a common problem with many GOF approaches is that the asymptotic distribution of the test statistic under the null hypothesis $C \in \mathcal{C}$ depends on the particular family $\mathcal{C}$ that is tested for, as well as on the unknown true value of the parameter $\boldsymbol{\lambda}$. In other words, many of the proposed GOF tests in the literature are not asymptotically \emph{distribution-free}. As a result, the asymptotic distribution of the test statistics under the null hypothesis cannot be tabulated for universal use, and approximate $p$-values have to be computed for each model via, e.g., specialized bootstrap procedures such as the ones outlined in \cite{genest:etal:2009}, App.\ A-D.

In this paper, we develop an approach to construct asymptotically distri\-bution-free GOF tests for any parametric copula family satisfying some rather mild smoothness assumptions. We do not propose a particular test statistic, but instead construct a whole test \emph{process} on the unit hypercube $[0,1]^d$, which converges weakly to a standard $d$-variate Wiener process under the null hypothesis. Thus GOF tests can be conducted by comparing the observed path of this test process with the statistical behavior of a standard Wiener process. Various functionals of the test process can be used for this comparison, such as the absolute maximum over $[0,1]^d$ or integral functionals. Since the weak limit of the test process is a standard process independent of the family $\mathcal{C}$ or the true value of $\boldsymbol{\lambda}$, the limiting distributions of these functionals will also be independent of $\mathcal{C}$ and $\boldsymbol{\lambda}$, and they only need to be tabulated once for use in all testing problems. In practice, this is an important advantage. Our results can also be used to test the GOF of fully specified copula models rather than parametric families.

We also show that our approach is optimal, in the sense that the obtained test process does not ``lose any information'' asymptotically. More precisely: when considering a sequence of contiguous alternatives approaching the null model, the distance in variation between the limiting processes under the null and the alternatives is as large as the limiting distance in variation of the data themselves. As a consequence, for a given sequence of contiguous alternatives, we can find a functional of our test process that yields an asymptotically optimal test and hence outperforms (or matches) competing procedures. Such a test retains the distribution-freeness advantage and thus avoids resampling procedures. Naturally, a variety of high power omnibus tests can be constructed as indicated above; see also Section \ref{SDA}.

Our approach relies on parametric estimation of the marginal distributions $F_1, \allowbreak \ldots, F_d$. That is, we assume that there is a parametric family of univariate dfs $\mathcal{F} = \{F_{\boldsymbol{\theta}} : \boldsymbol{\theta} \in \Theta\}$ such that $F_j \in \mathcal{F}$ for $j=1,\ldots,d$. In fact, this requirement can be relaxed to $F_j \in \mathcal{F}_j$ for $j=1,\ldots,d$, where the parametric families $\mathcal{F}_j$ may be different, but we will stick with $F_j \in \mathcal{F}$ for notational simplicity. The parametric structure of the margins might be naturally provided by the specifics of the data-generating process, or it might follow from theoretical considerations such as limit theorems, or it might be assumed based on independent empirical analysis or expert judgement. In this paper, we take the parametric structure of the margins as given, and focus on inference about the unknown copula.

The remainder of this paper is structured as follows. In Section 2, we introduce two estimators for the copula $C$, a parametric one that works under the null hypothesis $C \in \mathcal{C}$ and a semi-parametric one that works in general. We consider the normalized difference $\eta_n$ between these estimators and determine its weak limit $\eta$ under the null hypothesis, as the sample size $n$ tends to infinity. We will see that the distribution of $\eta$ depends on the family $\mathcal{C}$, as well as on the true values of the parameter $\boldsymbol{\lambda}$ and the marginal parameters, so $\eta_n$ cannot be used as a basis for distribution-free testing. The crucial step in our approach is introduced in Section 3, where we describe a transformation that turns $\eta$ into a standard Wiener process on $[0,1]^d$. In Section 4, we apply an empirical version of this transformation to $\eta_n$, and show that the resulting process $W_n$ converges weakly to a standard Wiener process on $[0,1]^d$ under the null hypothesis. This is our first main result, and the transformed process $W_n$ is the test process that was alluded to above. In Section 5, we investigate the behavior of the test process $W_n$ under a sequence of contiguous alternatives, and we show that transforming the raw data into $W_n$ does not lead to any loss of information asymptotically. This is our second main result. In Section 6, we present some simulation results that demonstrate the finite-sample behavior of some functionals of $W_n$ under the null and alternative hypotheses, and we then apply our approach to a real-world data set and analyze the results. Section 7 and an online appendix contain all the proofs.

\section{Comparing two copula estimators}\label{sec.2}
As in the Introduction, we assume that $\mathbf{X}_i = (X_{i1}, \ldots, X_{id})^\mathsf{T}, i \in \{1, \ldots, n\},$ are i.i.d.\ random vectors with common df $F$, which has continuous marginal dfs $F_1, \ldots, F_d$ and copula $C$. We further assume throughout that the marginal dfs are members of some parametric family of univariate dfs, $\mathcal{F} = \{ F_{\boldsymbol{\theta}} : \boldsymbol{\theta} \in \Theta\}$, indexed by $\boldsymbol{\theta} = (\theta_1, \ldots, \theta_m)^\mathsf{T} \in \Theta$, where $\Theta$ is some open subset of $\mathbb{R}^m$. This means that there exist $\boldsymbol{\theta}_1, \ldots, \boldsymbol{\theta}_d \in \Theta$ such that $F_j = F_{\boldsymbol{\theta}_j}$ for $j \in \{1, \ldots, d\}$. Our ultimate aim is to test the hypothesis $C \in \mathcal{C} = \{ C_{\boldsymbol{\lambda}} : \boldsymbol{\lambda} \in \Lambda\}$, where the parameter $\boldsymbol{\lambda} = (\lambda_1, \ldots, \lambda_p)^\mathsf{T}$ takes values in some open subset $\Lambda$ of $\mathbb{R}^p$. Throughout Sections 2-4, we will assume that the null hypothesis holds, i.e.\ there is a $\boldsymbol{\lambda}_0 \in \Lambda$ such that $C = C_{\boldsymbol{\lambda_0}}$.

There are various ways of estimating the copula from i.i.d.\ data, depending on the assumptions one is willing to make about the underlying model, as well as the requirements one chooses to impose on the estimator, such as smoothness. Perhaps the most straightforward and well-known copula estimator is the non-parametric \emph{empirical copula} discussed in \cite{ruymgaart:1973}, \cite{ruschendorf:1976}, \cite{deheuvels:1979, deheuvels:1981}, \cite{gaenssler:stute:1987}, \cite{fermanian:etal:2004} and \cite{segers:2012}, among others. Other commonly used approaches for copula estimation include two-step methods where the first step involves (non-parametric or parametric) estimation of the margins and the second step estimates the copula parametrically from marginal data transformed in accordance with the first step. Such estimators are studied in, e.g., \cite{genest:etal:1995} and \cite{shih:louis:1995}. A broad overview of various copula estimation methods can be found in \cite{charpentier:etal:2007}, \cite{choros:etal:2010} and Ch.\ 5 of \cite{joe:2015}.

In this paper, we will make use of two estimators for $C$: a parametric estimator $C_{\widehat{\boldsymbol{\lambda}}}$ and a semi-parametric estimator $\widehat{C}$. We do not specify the estimator $\widehat{\boldsymbol{\lambda}}$, but require it to satisfy a rather non-restrictive convergence assumption to be stated below. The semi-parametric estimator $\widehat{C}$ is defined as
\begin{equation*}
\widehat{C}(\mathbf{u}) = F_n(Q_{\widehat{\boldsymbol{\theta}}_1}(u_1), \ldots, Q_{\widehat{\boldsymbol{\theta}}_d}(u_d)), \quad \mathbf{u} \in [0,1]^d,
\end{equation*}
where $\widehat{\boldsymbol{\theta}}_1, \ldots, \widehat{\boldsymbol{\theta}}_d$ denote appropriate estimators for $\boldsymbol{\theta}_1, \ldots, \boldsymbol{\theta}_d$, $Q_{\boldsymbol{\theta}}$ denotes the quantile function associated with $F_{\boldsymbol{\theta}}$, and $F_n$ denotes the $d$-variate empirical df generated by the sample $\mathbf{X}_1, \ldots, \mathbf{X}_n$:
\begin{equation*}
F_n(\mathbf{x}) = \frac{1}{n} \sum_{i=1}^n \mathds{1} \{ \mathbf{X}_i \le \mathbf{x} \}, \quad \mathbf{x} \in \mathbb{R}^d.
\end{equation*}
Here, $\mathbf{X}_i \le \mathbf{x}$ is short-hand notation for ``$X_{ij} \le x_j$ for all $j=1,\ldots, d$''. Note that in view of the representation (\ref{cop2}) of the copula $C$ and our parametric assumption on the marginal dfs of $F$, the estimator $\widehat{C}$ is a natural one. To our knowledge, its asymptotic behavior has not been studied in the existing literature.

Under the null hypothesis, both $\widehat{C}$ and $C_{\widehat{\boldsymbol{\lambda}}}$ estimate the true copula $C$, while only $\widehat{C}$ correctly estimates $C$ when the null hypothesis does not hold. Thus the asymptotic discrepancy between the two estimators provides a natural starting point for a GOF test. With that in mind, we define
\begin{equation}\label{eta_n.def}
\eta_n(\mathbf{u}) = \sqrt{n} [\widehat{C}(\mathbf{u}) - C_{\widehat{\boldsymbol{\lambda}}}(\mathbf{u})], \quad \mathbf{u} \in [0,1]^d.
\end{equation}
Our first result will be a theorem describing the asymptotic behavior of $\eta_n$, but first we establish some notation and state the necessary assumptions about the various estimators introduced above, as well as about the parametric families $\mathcal{C}$ and $\mathcal{F}$.

Let $C_n$ denote the empirical df generated by the (unobserved) copula sample $(F_1(X_{i1}), \allowbreak \ldots, F_d(X_{id}))^\mathsf{T}$, $i \in \{1,\ldots,n\}$. That is,
\begin{equation*}
C_n(\mathbf{u}) = \frac{1}{n} \sum_{i=1}^n \mathds{1}\{F_1(X_{i1}) \le u_1, \ldots, F_d(X_{id}) \le u_d\}, \quad \mathbf{u} \in [0,1]^d.
\end{equation*}
Note that we can then write
\begin{equation*}
\widehat{C}(\mathbf{u}) = C_n\big(F_{\boldsymbol{\theta}_1}(Q_{\boldsymbol{\widehat{\theta}}_1}(u_1)), \ldots, F_{\boldsymbol{\theta}_d}(Q_{\boldsymbol{\widehat{\theta}}_d}(u_d))\big), \quad \mathbf{u} \in [0,1]^d.
\end{equation*}
We also define
\begin{equation}\label{alpha_n}
\alpha_n(\mathbf{u}) = \sqrt{n}[C_n(\mathbf{u}) - C(\mathbf{u})], \quad \mathbf{u} \in [0,1]^d,
\end{equation}
so that $\alpha_n$ is the classical empirical process associated with the df $C$. The asymptotic behavior of $\alpha_n$ is well-known, see e.g.\ \cite{neuhaus:1971}: we have $\alpha_n \Rightarrow B_C$ in the Skorohod space $D([0,1]^d)$, where ``$\Rightarrow$'' denotes weak convergence and $B_C$ is a $C$-Brownian bridge, that is, a mean-zero Gaussian process on $[0,1]^d$ with covariance structure
\begin{equation*}
E[B_C(\mathbf{u})B_C(\mathbf{u}')] = C(\mathbf{u} \wedge \mathbf{u}') - C(\mathbf{u})C(\mathbf{u}').
\end{equation*}
Here, $\mathbf{u} \wedge \mathbf{u}' := (u_1 \wedge u'_1, \ldots, u_d \wedge u'_d)^\mathsf{T}$.

The assumptions needed for our first result are listed below, followed by the result itself.\\

\textbf{A1.} There exist a $p$-variate random vector $\boldsymbol{\zeta}_0$ and $m$-variate random vectors $\boldsymbol{\zeta}_1, \ldots, \boldsymbol{\zeta}_d$ such that
\begin{equation}\label{joint.weak}
(\alpha_n, \sqrt{n}(\boldsymbol{\lambda}_0 - \widehat{\boldsymbol{\lambda}}), \sqrt{n}(\boldsymbol{\theta}_1 - \widehat{\boldsymbol{\theta}}_1), \ldots, \sqrt{n}(\boldsymbol{\theta}_d - \widehat{\boldsymbol{\theta}}_d)) \Rightarrow (B_C, \boldsymbol{\zeta}_0, \boldsymbol{\zeta}_1, \ldots, \boldsymbol{\zeta}_d)
\end{equation}
in $D([0,1]^d) \times \mathbb{R}^p \times (\mathbb{R}^m)^d$.\\

\textbf{A2.} The mappings
\begin{equation*}
(\mathbf{u},\boldsymbol{\lambda}) \mapsto \nabla C_{\boldsymbol{\lambda}}(\mathbf{u}) = \big( C_{\boldsymbol{\lambda}}^{(1)}(\mathbf{u}), \ldots, C_{\boldsymbol{\lambda}}^{(d)}(\mathbf{u}) \big)^\mathsf{T} := \Big( \frac{\partial C_{\boldsymbol{\lambda}}(\mathbf{u})}{\partial u_1}, \ldots, \frac{\partial C_{\boldsymbol{\lambda}}(\mathbf{u})}{\partial u_d} \Big)^\mathsf{T}
\end{equation*}
and
\begin{equation*}
(\mathbf{u},\boldsymbol{\lambda}) \mapsto \dt{C}_{\boldsymbol{\lambda}}(\mathbf{u}) = \big( \dt{C}_{\boldsymbol{\lambda}}^{(1)}(\mathbf{u}), \ldots, \dt{C}_{\boldsymbol{\lambda}}^{(p)}(\mathbf{u})\big)^\mathsf{T} := \Big( \frac{\partial C_{\boldsymbol{\lambda}}(\mathbf{u})}{\partial \lambda_1}, \ldots, \frac{\partial C_{\boldsymbol{\lambda}}(\mathbf{u})}{\partial \lambda_p} \Big)^\mathsf{T}
\end{equation*}
are continuous on $(0,1)^d \times \Lambda$.\\

\textbf{A3.} The mapping
\begin{equation*}
(x, \boldsymbol{\theta}) \mapsto \dt{F}_{\boldsymbol{\theta}}(x) = \big( \dt{F}_{\boldsymbol{\theta}}^{(1)}(x), \ldots, \dt{F}_{\boldsymbol{\theta}}^{(m)}(x) \big)^\mathsf{T} := \Big( \frac{\partial F_{\boldsymbol{\theta}}(x)}{\partial \theta_1}, \ldots, \frac{\partial F_{\boldsymbol{\theta}}(x)}{\partial \theta_m} \Big)^\mathsf{T}
\end{equation*}
is continuous on $\mathbb{R} \times \Theta$, the mapping $(u, \boldsymbol{\theta}) \mapsto Q_{\boldsymbol{\theta}}(u)$ is bounded on compact subsets of $(0,1) \times \Theta$, and the mapping $(u, \boldsymbol{\theta}) \mapsto \dt{F}_{\boldsymbol{\theta}}(Q_{\boldsymbol{\theta}}(u))$ is continuous on $(0,1) \times \Theta$.

\begin{thm}\label{thm.1}
Let $\eta_n$ be the process defined in (\ref{eta_n.def}), and let $0 < \delta < \tau < 1$. Under Assumptions A1-A3,
\begin{equation}\label{first.conv}
\begin{split}
\eta_n(\mathbf{u})
&\Rightarrow B_C(\mathbf{u}) + \sum_{j=1}^d C^{(j)}(\mathbf{u}) \dt{F}_{\boldsymbol{\theta}_j}(Q_{\boldsymbol{\theta}_j}(u_j))^\mathsf{T} \boldsymbol{\zeta}_j + \dt{C}(\mathbf{u})^\mathsf{T} \boldsymbol{\zeta}_0 \\
&=: \eta(\mathbf{u})
\end{split}
\end{equation}
in $D([\delta,\tau]^d)$, where $C^{(j)}$ and $\dt{C}$ are short-hand notation for $C^{(j)}_{\boldsymbol{\lambda}_0}$ and $\dt{C}_{\boldsymbol{\lambda}_0}$.
\end{thm}

\begin{rem}
Heuristically, the limiting process $\eta$ consists of a $C$-Brown\-ian bridge $B_C$ plus $dm$ additive terms ``contributed by'' the estimation of the $m$-dimensional parameters $\boldsymbol{\theta}_1, \ldots, \boldsymbol{\theta}_d$, plus $p$ additive terms ``contributed by'' the estimation of the $p$-dimensional parameter $\boldsymbol{\lambda}_0$. Since the distribution of $\eta$ clearly depends on the underlying family $\mathcal{C}$, as well as the (unknown) true values of $\boldsymbol{\theta}_1, \ldots, \boldsymbol{\theta}_d$ and $\boldsymbol{\lambda}_0$, Theorem \ref{thm.1} is far from suitable for distribution-free testing.
\end{rem}

\begin{rem}
The convergence in (\ref{first.conv}) does not necessarily hold in the space $D([0,1]^d)$ under the stated assumptions. For example, in the case $m=1$, consider the parametric family $\mathcal{F} = \{F_\theta : \theta \in (0,\infty)\}$, with
\begin{equation*}
F_\theta(x) = 1 - \sqrt{1-x/\theta}, \quad x \in (0,\theta).
\end{equation*}
Note that $F_\theta$ is a beta distribution with shape parameters fixed at 1 and $1/2$, and a free scale parameter $\theta > 0$. Also note that $F_\theta$ satisfies Assumption A3. However, for any $\theta > 0$, the expression
$$\dt{F}_\theta(Q_\theta(u)) = \frac{1}{2\theta}\Big( 1-u-\frac{1}{1-u} \Big)$$
is unbounded near $u=1$ and hence the process $\eta$ is in general not well-defined on the closed hypercube $[0,1]^d$, so the convergence in (\ref{first.conv}) cannot hold in $D([0,1]^d)$.
\end{rem}

\section{Transforming $\eta$ into a standard Wiener process}
As we observed in the previous section, the empirical process $\eta_n$ cannot directly be used as a basis for distribution-free testing, since its limiting process $\eta$ depends on the underlying family $\mathcal{C}$ and unknown parameter values. We will remedy this problem by transforming $\eta$ into a standard $d$-variate Wiener process. The transformation itself will depend on $\mathcal{C}$ and parameter values, but the distribution of the resulting process will not, which will facilitate asymptotically distribution-free testing. We first introduce some notation and assumptions.

Recall that
\begin{equation*}
B_C \stackrel{d}{=} V_C - CV_C(\mathbf{1}),
\end{equation*}
with $V_C$ a $C$-Wiener process on $[0,1]^d$, i.e.\ a mean-zero Gaussian process with covariance $E[V_C(\mathbf{u})V_C(\mathbf{u}')] \allowbreak = C(\mathbf{u} \wedge \mathbf{u}')$. We can thus alternatively express the limiting process $\eta$ in (\ref{first.conv}) as
\begin{equation}\label{eta}
\eta(\mathbf{u}) = V_C(\mathbf{u}) - C(\mathbf{u})V_C(\mathbf{1}) + \sum_{j=1}^d C^{(j)}(\mathbf{u}) \dt{F}_{\boldsymbol{\theta}_j} (Q_{\boldsymbol{\theta}_j}(u_j))^\mathsf{T} \boldsymbol{\zeta}_j + \dt{C}(\mathbf{u})^\mathsf{T} \boldsymbol{\zeta}_0,
\end{equation}
for $\mathbf{u} \in (0,1)^d$. Hence we see that $\eta$ is of the form
\begin{equation}\label{eta.2}
\eta(\mathbf{u}) = V_C(\mathbf{u}) + \sum_{i=1}^{1+dm+p} K_i(\mathbf{u}) Z_i, \quad \mathbf{u} \in (0,1)^d,
\end{equation}
where the $Z_i$ are some random variables, and the $K_i$ are deterministic functions on $(0,1)^d$ defined by
\begin{align}\label{K_i.def}
K_1(\mathbf{u}) &= C(\mathbf{u}), \notag\\
K_{1+(j-1)m+i}(\mathbf{u}) &= C^{(j)}(\mathbf{u}) \dt{F}^{(i)}_{\boldsymbol{\theta}_j}(Q_{\boldsymbol{\theta}_j}(u_j)), \quad j = 1, \ldots, d, \, i = 1, \ldots, m,\\
K_{1+dm+i}(\mathbf{u}) &= \dt{C}^{(i)}(\mathbf{u}), \quad i=1,\ldots,p. \notag
\end{align}
We note that (\ref{eta.2}) is analogous to the bivariate form (23) in \cite{can:etal:2015}. In that paper, a transformation of such processes into a standard bivariate Wiener process was described, which was an application of the ``innovation martingale transform'' idea developed in \cite{khmaladze:1981, khmaladze:1988, khmaladze:1993}. This idea has been applied to various statistical problems in the literature over the last couple of decades; see, for example, \cite{mckeague:etal:1995}, \cite{nikabadze:stute:1997}, \cite{stute:etal:1998}, \cite{koenker:xiao:2002, koenker:xiao:2006}, \cite{khmaladze:koul:2004, khmaladze:koul:2009}, \cite{delgado:etal:2005} and \cite{dette:hetzler:2009}. We will develop a suitable innovation martingale transform to construct a standard $d$-variate Wiener process on $[0,1]^d$ from the process $\eta$ in (\ref{eta.2}).

The approach here is novel in the sense that direct use of the data naturally leads to processes on the unit hypercube $[0,1]^d$, whereas in other applications of the martingale transform in multivariate contexts, the data are first \emph{transformed} to $[0,1]^d$ by an arbitrary transformation which influences the statistical properties of the procedures (see, for example, \cite{khmaladze:1993} and \cite{einmahl:khmaladze:2001}).

We first state the necessary assumptions and establish some notation.\\

\textbf{A4.} For each $\boldsymbol{\lambda} \in \Lambda$, the copula $C_{\boldsymbol{\lambda}}$ has a strictly positive density $c_{\boldsymbol{\lambda}}$ on $(0,1)^d$, and the mappings
\begin{equation*}
(\mathbf{u},\boldsymbol{\lambda}) \mapsto \nabla c_{\boldsymbol{\lambda}}(\mathbf{u}) = \big( c_{\boldsymbol{\lambda}}^{(1)}(\mathbf{u}), \ldots, c_{\boldsymbol{\lambda}}^{(d)}(\mathbf{u}) \big)^\mathsf{T} := \Big( \frac{\partial c_{\boldsymbol{\lambda}}(\mathbf{u})}{\partial u_1}, \ldots, \frac{\partial c_{\boldsymbol{\lambda}}(\mathbf{u})}{\partial u_d} \Big)^\mathsf{T}
\end{equation*}
and
\begin{equation*}
(\mathbf{u},\boldsymbol{\lambda}) \mapsto \dt{c}_{\boldsymbol{\lambda}}(\mathbf{u}) = \big( \dt{c}_{\boldsymbol{\lambda}}^{(1)}(\mathbf{u}), \ldots, \dt{c}_{\boldsymbol{\lambda}}^{(p)}(\mathbf{u})\big)^\mathsf{T} := \Big( \frac{\partial c_{\boldsymbol{\lambda}}(\mathbf{u})}{\partial \lambda_1}, \ldots, \frac{\partial c_{\boldsymbol{\lambda}}(\mathbf{u})}{\partial \lambda_p} \Big)^\mathsf{T}
\end{equation*}
are continuous on $(0,1)^d \times \Lambda$.\\

Now, with the functions $K_i$ as defined in (\ref{K_i.def}), let us denote
\begin{equation}\label{k_i.def}
k_i(\mathbf{u}) = \dd K_i(\mathbf{u})/\dd C(\mathbf{u}), \quad i=1,\ldots,1+dm+p,
\end{equation}
so that
\begin{align*}
k_1(\mathbf{u})
&= 1,\\
k_{1+(j-1)m+i}(\mathbf{u})
&= \dt{F}^{(i)}_{\boldsymbol{\theta}_j}(Q_{\boldsymbol{\theta}_j}(u_j)) \frac{\partial}{\partial u_j} \log c_{\boldsymbol{\lambda}_0}(\mathbf{u}) + \frac{\partial}{\partial u_j} \dt{F}^{(i)}_{\boldsymbol{\theta}_j}(Q_{\boldsymbol{\theta}_j}(u_j)),\\
& \hspace{125pt} j = 1, \ldots, d, \, i = 1, \ldots, m,\\
k_{1+dm+i}(\mathbf{u})
&= \frac{\partial}{\partial \lambda_i} \log c_{\boldsymbol{\lambda}}(\mathbf{u}) \Big|_{\boldsymbol{\lambda} = \boldsymbol{\lambda}_0}, \quad i=1,\ldots,p.
\end{align*}
Let $\mathbf{k}(\mathbf{u})$ denote the column vector consisting of $k_1(\mathbf{u}),\ldots,k_{1+dm+p}(\mathbf{u})$. We will also write $\mathbf{k}(\mathbf{u}, \boldsymbol{\theta}'_1, \ldots, \boldsymbol{\theta}'_d, \boldsymbol{\lambda}')$ for the vector $\mathbf{k}(\mathbf{u})$ with true parameter values $\boldsymbol{\theta}_1, \ldots, \boldsymbol{\theta}_d$, $\boldsymbol{\lambda}_0$ replaced by arbitrary values $\boldsymbol{\theta}'_1, \ldots, \boldsymbol{\theta}'_d \in \Theta, \boldsymbol{\lambda}' \in \Lambda$.

Finally, given $0 < \delta < 1/2$, let
\begin{equation*}
S_\delta(t) = [\delta, 1 - \delta/2]^{d-1} \times [t, 1-\delta/2], \quad t \in [\delta, 1-\delta/2),
\end{equation*}
and introduce matrices
\begin{equation}\label{I.def}
\mathbf{I}_\delta(t) = \int_{S_\delta(t)} \mathbf{k}(\mathbf{s}) \mathbf{k}(\mathbf{s})^\mathsf{T} \, \dd C(\mathbf{s}), \quad t \in [\delta, 1-\delta/2).
\end{equation}
We also define
\begin{equation}\label{I.def.2}
\begin{split}
\mathbf{I}_\delta(t, \boldsymbol{\theta}'_1, \ldots, \boldsymbol{\theta}'_d, \boldsymbol{\lambda}') = \int_{S_\delta(t)} \mathbf{k}(\mathbf{s},\boldsymbol{\theta}'_1, \ldots, \boldsymbol{\theta}'_d, \boldsymbol{\lambda}') \, \mathbf{k}(\mathbf{s},\boldsymbol{\theta}'_1, \ldots, \boldsymbol{\theta}'_d, \boldsymbol{\lambda}')^\mathsf{T} \, \dd C_{\boldsymbol{\lambda}'}(\mathbf{s}),\\
 t \in [\delta, 1-\delta/2).
\end{split}
\end{equation}

Note that in the nomenclature of likelihood theory, $\mathbf{k}$ is the vector of \emph{score functions} for the underlying copula model (extended by the constant function 1 in the first component), and $\mathbf{I}_\delta(t)$ is a partial \emph{Fisher information matrix} constructed from these functions.

Our next assumption is:\\

\textbf{A5.} The matrices $\mathbf{I}_\delta(t, \boldsymbol{\theta}'_1, \ldots, \boldsymbol{\theta}'_d, \boldsymbol{\lambda}')$ in (\ref{I.def.2}) are well-defined and invertible for all $0 < \delta < 1/2$, $t \in [\delta, 1-\delta/2)$, $\boldsymbol{\theta}'_1, \ldots, \boldsymbol{\theta}'_d \in \Theta, \boldsymbol{\lambda}' \in \Lambda$.\\

Given $0 < \delta < 1/2$, let $\boldsymbol{\delta}$ denote the point $(\delta, \ldots, \delta)^\mathsf{T} \in (0,1)^d$, and given points $\mathbf{u}, \mathbf{v} \in [0,1]^d$, let $[\mathbf{u},\mathbf{v}]$ denote the hyperrectangle $[u_1,v_1] \times \ldots \times [u_d,v_d]$. Also, given $a > 0$, let $a\mathbf{u}$ denote the point $(au_1, \ldots, au_d)^\mathsf{T}$.

We are now ready to state our transformation result.
\begin{thm}\label{thm.2}
Let $\eta$ be the limiting process appearing in Theorem \ref{thm.1} and let $0 < \delta < 1/2$. If Assumptions A4-A5, restricted to the true parameter values $\boldsymbol{\theta}_1, \ldots, \boldsymbol{\theta}_d$, $\boldsymbol{\lambda}_0$, hold, then the process
\begin{align}\label{W.def}
W(\mathbf{u}) &= \frac{1}{(1-2\delta)^{d/2}} \Bigg[ \int_{[\boldsymbol{\delta}, \boldsymbol{\delta}+(1-2\delta)\mathbf{u}]} \frac{1}{\sqrt{c(\mathbf{s})}} \, \dd \eta(\mathbf{s}) \notag \\
& {} \hspace{40pt} - \int_{[\boldsymbol{\delta}, \boldsymbol{\delta}+(1-2\delta)\mathbf{u}]} \mathbf{k}(\mathbf{s})^\mathsf{T}
\bigg( \mathbf{I}_\delta^{-1}(s_d) \int_{S_\delta(s_d)} \mathbf{k}(\mathbf{s'})\,\dd \eta(\mathbf{s'}) \bigg) \sqrt{c(\mathbf{s})} \, \dd \mathbf{s} \Bigg]
\end{align}
is a standard Wiener process on $[0,1]^d$.
\end{thm}

\begin{rem}
The transformation in (\ref{W.def}) ``annihilates'' the terms following $V_C$ in (\ref{eta.2}) to produce a $C$-Wiener process on $[\delta, 1-\delta]^d$, which is then normalized and scaled to the entire hypercube $[0,1]^d$, so that the end result is indeed a standard Wiener process on $[0,1]^d$. In the next section, we will describe how the transformation of Theorem \ref{thm.2} facilitates asymptotically distribution-free testing for $C$.
\end{rem}

\section{Goodness-of-fit testing: null hypothesis}\label{simp.hyp}

Recall the empirical process $\eta_n$ defined in (\ref{eta_n.def}). In Theorem \ref{thm.1} we have derived its weak limit $\eta$ as $n \to \infty$, and in Theorem \ref{thm.2} we have described a transformation that turns $\eta$ into a standard Wiener process on $[0,1]^d$. In this section, we will apply the same transformation (or rather, its empirical version, with unknown parameters replaced by estimators) to $\eta_n$, and we will show that the resulting process converges weakly to a standard Wiener process. This is the first main result of this paper.

Applying transformation (\ref{W.def}) to $\eta_n$, with unknown parameters replaced by estimators, we obtain the following empirical process on $[0,1]^d$:
\begin{align}\label{W_n.def}
W_n(\mathbf{u}) &= \frac{1}{(1-2\delta)^{d/2}} \Bigg[ \int_{[\boldsymbol{\delta}, \boldsymbol{\delta}+(1-2\delta)\mathbf{u}]} \frac{1}{\sqrt{c_{\widehat{\boldsymbol{\lambda}}}(\mathbf{s})}} \, \dd \eta_n(\mathbf{s}) \notag \\
& {} \hspace{23pt} - \int_{[\boldsymbol{\delta}, \boldsymbol{\delta}+(1-2\delta)\mathbf{u}]} \widehat{\mathbf{k}}(\mathbf{s})^\mathsf{T}
\bigg( \widehat{\mathbf{I}}_\delta^{-1}(s_d) \int_{S_\delta(s_d)} \widehat{\mathbf{k}}(\mathbf{s'})\,\dd \eta_n(\mathbf{s'}) \bigg) \sqrt{c_{\widehat{\boldsymbol{\lambda}}}(\mathbf{s})} \, \dd \mathbf{s} \Bigg].
\end{align}
Here, $\widehat{\mathbf{k}}(\cdot)$ and $\widehat{\mathbf{I}}_\delta(\cdot)$ are short-hand notations for $\mathbf{k}(\cdot, \widehat{\boldsymbol{\theta}}_1, \ldots, \widehat{\boldsymbol{\theta}}_d, \widehat{\boldsymbol{\lambda}})$ and $\mathbf{I}_\delta(\cdot,$ $\widehat{\boldsymbol{\theta}}_1, \ldots, \widehat{\boldsymbol{\theta}}_d, \widehat{\boldsymbol{\lambda}})$, respectively.

Before stating the convergence result on $W_n$, we introduce some further notation and assumptions. Given a hyperrectangle $[\mathbf{a}, \mathbf{b}] \subset \mathbb{R}^d$ and a function $\varphi: [\mathbf{a}, \mathbf{b}] \to \mathbb{R}$, let $V_{[\mathbf{a}, \mathbf{b}]}(\varphi)$ denote the total variation of $\varphi$ on $[\mathbf{a}, \mathbf{b}]$ in the sense of Vitali; see e.g.\ \cite{owen:2005}, Sec.\ 4 for a definition. Also, given $I \subset \{1,\ldots,d\}$ and $\mathbf{x} \in \mathbb{R}^d$, let $|I|$ denote the cardinality of $I$, and let $\mathbf{x}_I$ denote the point in $\mathbb{R}^{|I|}$ obtained by discarding all coordinates $x_j$ of $\mathbf{x}$ for $j \notin I$. Moreover, given disjoint subsets $I_1, I_2, I_3 \subset \{1,\ldots,d\}$ with $I_1 \cup I_2 \cup I_3 = \{1,\ldots,d\}$, let $\varphi(\mathbf{x}_{I_1};\mathbf{a}_{I_2},\mathbf{b}_{I_3})$ denote the function on $[\mathbf{a}_{I_1}, \mathbf{b}_{I_1}]$ obtained by fixing the $j^\text{th}$ argument of $\varphi$ at $a_j$ for $j \in I_2$ and at $b_j$ for $j \in I_3$.

We consider an alternative concept of ``total variation'' on $[\mathbf{a}, \mathbf{b}]$, as follows:
\begin{equation}\label{hk.def}
V^{\text{HK}}_{[\mathbf{a}, \mathbf{b}]}(\varphi) := \sum_{\substack{I_1, I_2, I_3 \subset \{1,\ldots,d\}, \, I_1 \neq \varnothing \\ I_1 + I_2 + I_3 = \{1,\ldots,d\}}} V_{[\mathbf{a}, \mathbf{b}]}(\varphi(\mathbf{x}_{I_1};\mathbf{a}_{I_2},\mathbf{b}_{I_3})),
\end{equation}
with $I_1 + I_2 + I_3$ denoting a disjoint union. In other words, $V^{\text{HK}}_{[\mathbf{a}, \mathbf{b}]}(\cdot)$ sums the Vitali variations over the hyperrectangle $[\mathbf{a}, \mathbf{b}]$ and over all of its ``faces'' where the $j^\text{th}$ coordinate is fixed at $a_j$ or $b_j$, for at least one $j \in \{1,\ldots,d\}$. Note that $V^{\text{HK}}_{[\mathbf{a}, \mathbf{b}]}$ is a variant of the so-called \emph{Hardy-Krause variation} for multivariate functions; cf.\ \cite{owen:2005}, Def.\ 2. For $0 < \delta < 1/2$ and functions $\varphi: [\delta, 1-\delta]^d \to \mathbb{R}$, we will write $V^{\text{HK}}_\delta$ instead of $V^{\text{HK}}_{[\delta,1-\delta]^d}$, for brevity.

Let us also denote
\begin{equation*}
\gamma(\mathbf{u}) = \frac{1}{\sqrt{c(\mathbf{u})}}, \quad \widehat{\gamma}(\mathbf{u}) = \frac{1}{\sqrt{c_{\widehat{\boldsymbol{\lambda}}}(\mathbf{u})}}, \quad \Delta\gamma(\mathbf{u}) = \widehat{\gamma}(\mathbf{u}) - \gamma(\mathbf{u}), \quad \mathbf{u} \in (0,1)^d.
\end{equation*}
Similarly, we will denote $\Delta k_i = \widehat{k}_i - k_i$ for $i = 1, \ldots, 1+dm+p$, with the $k_i$ as defined in (\ref{k_i.def}). We introduce the final assumption needed for our convergence result:\\

\textbf{A6.} For any $0 < \delta < 1/2$, we have $V^{\text{HK}}_\delta(\gamma) < \infty$ and $V^{\text{HK}}_\delta(\Delta\gamma) = o_P(1)$. Also, $V^{\text{HK}}_\delta(k_i) < \infty$ and $V^{\text{HK}}_\delta(\Delta k_i) = o_P(1)$ for $i = 1, \ldots, 1+dm+p$.\\

\begin{thm}\label{thm.3}
Let $0 < \delta < 1/2$. Under Assumptions A1-A6, the process $W_n$ in (\ref{W_n.def}) converges weakly to a standard Wiener process in $D([0,1]^d)$.
\end{thm}

\begin{rem}
Theorem \ref{thm.3} is analogous to Theorem \ref{thm.1} in that it describes the asymptotic behavior of an empirical process constructed from the data $\mathbf{X}_1, \ldots, \mathbf{X}_n$. However, unlike the process $\eta_n$, the asymptotic behavior of $W_n$ is distribution-free: it converges to a standard Wiener process. Thus a test for the null hypothesis can now be performed by assessing how the observed path of $W_n$ compares to the ``usual'' statistical behavior of a standard Wiener process. Since this comparison can be done through many different functionals of $W_n$, we can construct a multitude of asymptotically distribution-free tests. In Section \ref{SDA}, we demonstrate through simulations and a real-world data analysis how such tests can be conducted.
\end{rem}

\begin{rem}
The statement and proof of Theorem \ref{thm.3} also applies, with obvious modifications, in the case $\mathcal{C}=\{C_0\}$, where $C_0$ denotes a fully specified copula. Thus the test process $W_n$ can also be used for testing null hypotheses of the form $C=C_0$. This will also be demonstrated in the simulations of Section \ref{SDA}.
\end{rem}

\section{Goodness-of-fit testing: contiguous alternatives}

We now consider testing $C \in \mathcal{C} = \{C_{\boldsymbol{\lambda}}: \boldsymbol{\lambda} \in \Lambda\}$ when the true copula of the underlying sample does not lie in $\mathcal{C}$ but approaches it as the sample size grows.

So let us assume that, for each $n \ge 1$, we have an i.i.d.\ sample
\begin{equation}\label{n.samp}
\mathbf{X}_{(n)1} = (X_{(n)11}, \ldots, X_{(n)1d})^\mathsf{T}, \ldots, \mathbf{X}_{(n)n} = (X_{(n)n1}, \ldots, X_{(n)nd})^\mathsf{T}
\end{equation}
generated from a $d$-variate df $F_{(n)}$ with continuous margins $F_1, \ldots, F_d$ and copula $C_{(n)}$. We assume that the marginal dfs are independent of $n$, and (as in the previous sections) they are all members of some parametric family $\mathcal{F} = \{F_{\boldsymbol{\theta}}: \boldsymbol{\theta} \in \Theta\}$, so that there are $\boldsymbol{\theta}_1, \ldots, \boldsymbol{\theta}_d \in \Theta$ with $F_j = F_{\boldsymbol{\theta}_j}$ for $j=1,\ldots,d$. Regarding the sequence of copulas $C_{(1)}, C_{(2)}, \ldots$, we assume the following:\\

\textbf{B0.} There exists $\boldsymbol{\lambda}_0 \in \Lambda$ such that
\begin{equation*}
\bigg[\frac{\dd C_{(n)}}{\dd C_{\boldsymbol{\lambda}_0}}\bigg]^{1/2} = 1 + \frac{1}{2\sqrt{n}} h_n, \quad n=1,2,\ldots,
\end{equation*}
for a sequence of functions $h_1, h_2, \ldots$ supported on $[\delta,1-\delta]^d$, for some $0 < \delta < 1/2$. The functions $h_n$ satisfy
\begin{equation*}
\int_{[0,1]^d} (h_n-h)^2 \dd C_{\boldsymbol{\lambda}_0} \to 0 \text{ as } n \to \infty,
\end{equation*}
for some function $h$ with
\begin{equation*}
\int_{[0,1]^d} h^2 \, \dd C_{\boldsymbol{\lambda}_0} \in (0, \infty), \quad \int_{[0,1]^d} k_i h \, \dd C_{\boldsymbol{\lambda}_0} = 0,
\end{equation*}
where the functions $k_i, i=1,\ldots,1+dm+p,$ are as defined in (\ref{k_i.def}).\\

Note that for each $n \ge 1$, the distribution of the sample in (\ref{n.samp}) on $(\mathbb{R}^d)^n$ is given by the $n$-fold product measure $F^n_{(n)} = F_{(n)} \times \ldots \times F_{(n)}$, whereas if the underlying copula was equal to $C_{\boldsymbol{\lambda}_0}$, this distribution would of course be $F^n_0 =  F_0 \times \ldots \times F_0$, with $F_0$ denoting the df with margins $F_1, \ldots, F_d$ and copula $C_{\boldsymbol{\lambda}_0}$. It follows from \cite{oosterhoff:vanzwet:1979} that condition B0 is sufficient to make the sequence $\{F^n_{(n)}\}$ \emph{contiguous} with respect to $\{F^n_0\}$, in the sense that $\lim_{n \to \infty} F^n_0(A_n) = 0$ implies $\lim_{n \to \infty} F^n_{(n)}(A_n) = 0$, for any sequence of measurable sets $A_n \subset (\mathbb{R}^d)^n$.

Our first result in this section will establish the asymptotic behavior of $\eta_n$ in (\ref{eta_n.def}) in the present setting. We define, analogously to (\ref{alpha_n}),
\begin{equation*}
\alpha_{n}(\mathbf{u}) = \sqrt{n}[C_n(\mathbf{u}) - C_{(n)}(\mathbf{u})], \quad \mathbf{u} \in [0,1]^d,
\end{equation*}
where $C_n$ is the empirical df generated by the (unobserved) copula sample
\begin{equation*}
\big(F_1(X_{(n)11}), \ldots, F_d(X_{(n)1d})\big)^\mathsf{T}, \ldots, \big(F_1(X_{(n)n1}), \ldots, F_d(X_{(n)nd}) \big)^\mathsf{T},
\end{equation*}
and we state the following analogue of Assumption A1 in Section \ref{sec.2}:\\

\textbf{B1.} There exist a $p$-variate random vector $\boldsymbol{\zeta}_a$ and $m$-variate random vectors $\boldsymbol{\zeta}_1, \ldots, \boldsymbol{\zeta}_d$ such that
\begin{equation*}
\begin{split}
\lefteqn{(\alpha_n, \sqrt{n}(\boldsymbol{\lambda}_0 - \widehat{\boldsymbol{\lambda}}), \sqrt{n}(\boldsymbol{\theta}_1 - \widehat{\boldsymbol{\theta}}_1), \ldots, \sqrt{n}(\boldsymbol{\theta}_d - \widehat{\boldsymbol{\theta}}_d))} \hspace{190pt}\\
&\Rightarrow (B_{C_{\boldsymbol{\lambda}_0}}, \boldsymbol{\zeta}_a, \boldsymbol{\zeta}_1, \ldots, \boldsymbol{\zeta}_d)
\end{split}
\end{equation*}
in $D([0,1]^d) \times \mathbb{R}^p \times (\mathbb{R}^m)^d$.\\

The following result is the analogue of Theorem \ref{thm.1}, and its proof, which we omit, follows along similar lines.

\begin{thm}\label{thm.1.a}
Let $\eta_n$ be the process defined in (\ref{eta_n.def}), and let $0 < \varepsilon < \tau < 1$. Under Assumptions B0-B1 and A2-A3,
\begin{equation*}
\begin{split}
\eta_n(\mathbf{u})
&\Rightarrow B_{C_{\boldsymbol{\lambda}_0}}(\mathbf{u}) + \sum_{j=1}^d C_{\boldsymbol{\lambda}_0}^{(j)}(\mathbf{u}) \dt{F}_{\boldsymbol{\theta}_j}(Q_{\boldsymbol{\theta}_j}(u_j))^\mathsf{T} \boldsymbol{\zeta}_j + \dt{C}_{\boldsymbol{\lambda}_0}(\mathbf{u})^\mathsf{T} \boldsymbol{\zeta}_a \\
&\hspace{190pt} + \int_{[\mathbf{0}, \mathbf{u}]} h(\mathbf{s}) \, \dd C_{\boldsymbol{\lambda}_0}(\mathbf{s})\\
&=: \eta_a(\mathbf{u}) + \int_{[\mathbf{0}, \mathbf{u}]} h(\mathbf{s}) \, \dd C_{\boldsymbol{\lambda}_0}(\mathbf{s})
\end{split}
\end{equation*}
in $D([\varepsilon,\tau]^d)$.
\end{thm}

Next, we establish the asymptotic behavior of the test process $W_n$ in (\ref{W_n.def}) in the present setting. If we let $W_a$ denote the process $W$ in (\ref{W.def}), with $\eta$ replaced by $\eta_a$, then $W_a$ is still a standard Wiener process on $[0,1]^d$, since the change of $\eta$ to $\eta_a$ does not affect the proof of Theorem \ref{thm.2}. This, together with Theorem \ref{thm.1.a} above, yields the following analogue of Theorem \ref{thm.3}, which shows that under the sequence of contiguous alternatives, $W_n$ converges to a standard Wiener process plus a deterministic shift term. 

\begin{thm}\label{thm.3.a}
Under Assumptions B0-B1 and A2-A6, the process $W_n$ in (\ref{W_n.def}) converges weakly to $\widetilde{W} := W+S$ in $D([0,1]^d)$, where
\begin{equation*}
S(\mathbf{u}) = \frac{1}{(1-2\delta)^{d/2}} \int_{[\boldsymbol{\delta}, \boldsymbol{\delta} + (1-2\delta)\mathbf{u}]} g(\mathbf{s}) \sqrt{c_{\boldsymbol{\lambda}_0}(\mathbf{s})} \, \dd \mathbf{s},
\end{equation*}
with
\begin{equation}\label{g.def}
g(\mathbf{s}) = h(\mathbf{s}) - \mathbf{k}(\mathbf{s})^\mathsf{T} \mathbf{I}_\delta^{-1}(s_d) \int_{S_\delta(s_d)} \mathbf{k}(\mathbf{s'})h(\mathbf{s'})\,\dd C_{\boldsymbol{\lambda}_0}(\mathbf{s'}), \quad \mathbf{s} \in [\delta, 1-\delta]^d.
\end{equation}
\end{thm}

In order to judge how ``sensitive'' the test process $W_n$ is to the sequence of alternatives $F_{(1)}, F_{(2)}, \ldots$, we first recall the notion of \emph{distance in variation} for probability measures. Given two such measures $P$ and $\widetilde{P}$ defined on some sigma-algebra $\mathcal{B}$, the distance in variation between $P$ and $\widetilde{P}$ is defined as
\begin{equation*}
d(P,\widetilde{P}) = \sup_{B \in \mathcal{B}} |P(B)-\widetilde{P}(B)|.
\end{equation*}
If $L_n$ denotes the log-likelihood ratio $\log (\dd F^n_{(n)}/\dd F_0^n)$, then we know from likelihood theory that
\begin{equation}\label{add.1}
d(F^n_{(n)},F_0^n) = F^n_{(n)}(L_n>0) - F_0^n(L_n>0).
\end{equation}
Moreover, we also know that as $n \to \infty$,
\begin{equation}\label{add.2}
L_n \stackrel{d}{\rightarrow} 
\begin{cases} 
N\big(-\frac{1}{2} \|h\|^2, \|h\|^2\big) &\text{ under } F_0^n,\\[5pt]
N\big(\frac{1}{2} \|h\|^2, \|h\|^2\big) &\text{ under } F^n_{(n)},
\end{cases}
\end{equation}
where $N(\mu,\sigma^2)$ denotes the normal distribution with mean $\mu$ and variance $\sigma^2$, and
\begin{equation}\label{add.3}
\|h\| := \bigg( \int_{[0,1]^d} h^2 \, \dd C_{\boldsymbol{\lambda}_0} \bigg)^{1/2}.
\end{equation}
Combining (\ref{add.1}) and (\ref{add.2}), we see that 
\begin{equation*}
d(F^n_{(n)}, F^n_0) \to \nu(h) := 2\Phi\Big(\frac{1}{2}\|h\|\Big)-1
\end{equation*}
as $n \to \infty$, with $\Phi$ denoting the standard normal cdf.

The following result is the second main result of the paper. It establishes asymptotic optimality of the test process $W_n$.

\begin{thm}\label{thm.4}
Let $Q$ denote the distribution of a standard Wiener process $W$ on $D([0,1]^d)$, and let $\widetilde{Q}$ denote the distribution of the process $\widetilde{W}$ defined in Theorem \ref{thm.3.a}. Then, 
\begin{equation*}
\log\Big(\frac{\dd\widetilde{Q}}{\dd Q}\Big) \sim 
\begin{cases} 
N\big(-\frac{1}{2} \|h\|^2, \|h\|^2\big) &\text{ under } Q,\\[5pt]
N\big(\frac{1}{2} \|h\|^2, \|h\|^2\big) &\text{ under } \widetilde{Q}.
\end{cases}
\end{equation*}
Hence $d(\widetilde{Q},Q) = \nu(h)$.
\end{thm}

\begin{rem}
The result shows that the limiting distance in variation of the processes $W_n$ under the null and the contiguous alternatives is the same as that of the samples: $\nu(h)$. In fact, the respective distributions of the log-likelihood ratio $\log(\dd\widetilde{Q}/\dd Q)$ under the two measures are identical to the limiting distributions of $\log(\dd F_{(n)}^n/\dd F_0^n)$. Hence, the process $W_n$ is asymptotically as good as the data themselves for testing purposes. Indeed, for a given sequence of alternatives satisfying Assumption B0, consider the test that rejects $H_0$ if 
\begin{align*}
\lefteqn{(1-2\delta)^{d/2}\int_{[0,1]^d} \widehat g\big(\boldsymbol{\delta} + (1-2\delta)\mathbf{s}\big)\sqrt{c_{\widehat{\boldsymbol{\lambda}}}\big(\boldsymbol{\delta} + (1-2\delta)\mathbf{s}\big)} \, \dd W_n(\mathbf{s})} \hspace{230pt}\\
&\ge \widehat{\|h\|}\Phi^{-1}(1-\alpha),
\end{align*}
with $\widehat{g}$ and $\widehat{\|h\|}$ denoting the obvious estimators of $g$ in (\ref{g.def}) and $\|h\|$ in (\ref{add.3}). Then under regularity assumptions the probability of a type I error converges to $\alpha$ as $n\to \infty$, and the power converges to $1-\Phi(\Phi^{-1}(1-\alpha)-||h||)$.
According to the Neyman-Pearson Lemma, this limiting power is equal to that of the most powerful level-$\alpha$ tests for a simple null (picked from our $H_0$) against the simple $h_n$-alternatives. Hence, for the more general problem of testing the composite null hypothesis $C \in \mathcal{C}$ against the composite $h_n$-alternatives, we have an asymptotically uniformly most powerful test. This optimality shows that our approach can favorably compete in terms of power with any other approach in the literature.
\end{rem}

\section{Simulations and data analysis}\label{SDA}

In this section we present the results of a simulation study and a data analysis in order to illustrate the applicability of our approach in finite samples. All computations are performed in \texttt{R}. The code to implement the simulations and the data analysis is available from the authors upon request. 

\subsection{Simulation study}\label{sim.stu}

We consider two widely used parametric copula models, namely Clayton and Gumbel, and we demonstrate how one might test for the goodness-of-fit of these models, both in parametric and fully specified form, using our approach. We limit the simulations to the bivariate case, and we perform tests on simulated data both under the null and alternative hypotheses.

To be more specific, we consider the following copula models:

\begin{itemize}
\setlength\itemsep{1em}
\item \emph{Clayton:} $C_\lambda(u,v) = (u^{-\lambda} + v^{-\lambda} -1)^{-1/\lambda}, \; \lambda \in (0,\infty)$
\item \emph{Clayton(2):} $C(u,v) = (u^{-2} + v^{-2} - 1)^{-1/2}$
\item \emph{Gumbel:} $C_\lambda(u,v) = \exp\{-[(-\log u)^\lambda + (-\log v)^\lambda]^{1/\lambda}\}, \; \lambda \in [1,\infty)$
\item \emph{Gumbel(2):} $C(u,v) = \exp\{-[(-\log u)^2 + (-\log v)^2]^{1/2}\}$
\end{itemize}

To test for Clayton and Clayton(2) models under the null hypothesis, we generate 1000 samples of size $n=200$ from the bivariate distribution with Exponential(1) margins and Clayton(2) copula. To test for Gumbel and Gumbel(2) models under the null hypothesis, we generate 1000 samples of size $n=200$ from the bivariate distribution with Lomax(3,1) margins and Gumbel(2) copula. Recall that for $\alpha>0$ and $\sigma>0$, the Lomax($\alpha$, $\sigma$) distribution has the cdf
$$
F(x) = 1 - \Big(1+\frac{x}{\sigma}\Big)^{-\alpha}, \quad x>0,
$$
so it is a shifted Pareto distribution with tail parameter $\alpha$ and scale parameter $\sigma$.

From each simulated sample, we compute the test process $W_n$ in (\ref{W_n.def}) on a $100 \times 100$ grid $\mathcal{G}$ of equally spaced points covering $(0,1)^2$. Parameter estimates are computed through maximum likelihood (ML) estimation. For the parametric Clayton and Gumbel models, ML estimates are computed for the entire bivariate distribution, while for the Clayton(2) and Gumbel(2) models, ML estimates are computed separately for the two marginal parameter sets. Note that Assumption A1 holds for these estimators, which follows from arguments similar to those for asymptotic normality of MLEs. Assumptions A2-A6 are smoothness assumptions that are straightforward to verify for the considered models.

To compare the observed paths of $W_n$ to a standard Wiener process, two functionals are computed from each path of $W_n$, namely:
\begin{align*}
\kappa_n &= \max_{(x,y) \in \mathcal{G}} \big|W_n(x,y)\big|, \hspace{32pt} \text{(Kolmogorov-Smirnov type statistic)}\\
\omega_n^2 &=\|\mathcal{G}\|^2 \sum_{(x,y) \in \mathcal{G}}  W_n(x,y)^2, \hspace{12pt} \text{(Cram\'{e}r-von Mises type statistic)}
\end{align*}
where $\|\mathcal{G}\|$ denotes the mesh length of the grid $\mathcal{G}$, i.e.~1/100. To create benchmark distribution tables for these statistics, we also simulate 10{,}000 true standard Wiener process paths on the grid $\mathcal{G}$ and compute the same functionals for each path. We denote these functionals by $\kappa$ and $\omega^2$.

\begin{figure}[t]
\epsfig{file=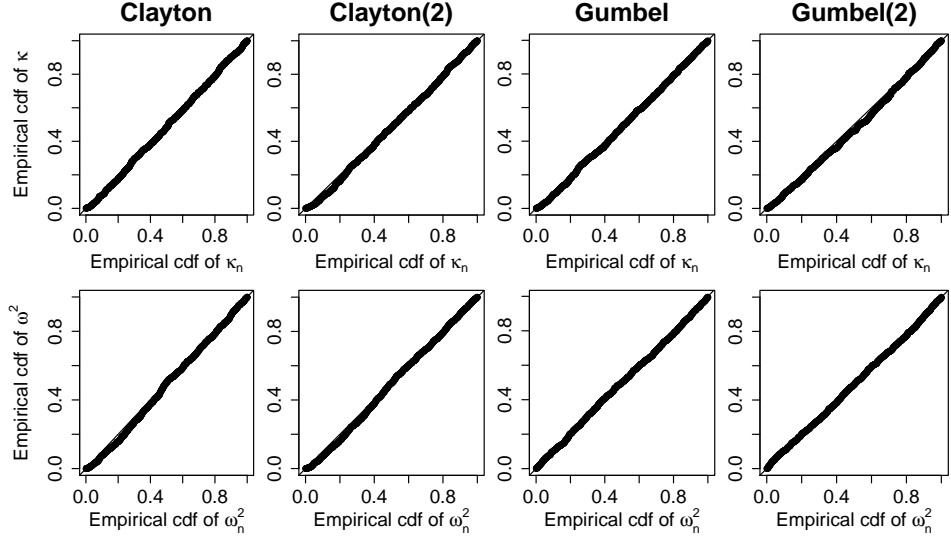, width=\textwidth}
\caption[]{PP-plots for the Kolmogorov-Smirnov (top) and Cram\'{e}r-von Mises (bottom) type test statistics}
\label{fig.qq}
\end{figure}

For each model, we construct PP-plots to compare the empirical distributions of $\kappa_n$ and $\omega^2_n$ with the theoretical distributions of $\kappa$ and $\omega^2$ (as inferred from the 10{,}000 simulated Wiener process paths). The results are shown in Fig.\ \ref{fig.qq}. We observe a very good match of empirical and limiting distributions for both statistics, especially in the upper right corners of the plots, which are important for testing. These results suggest that Theorem \ref{thm.3} yields good finite-sample approximations. This is confirmed by the observed fractions of Type I errors at 5\% and 1\% significance levels, given in Table \ref{size}. Note that the rejection counts are consistent with draws from a Binomial$(1000,0.05)$ or a Binomial$(1000,0.01)$ distribution, respectively.

We emphasize here that due to the distribution-free nature of our approach, the critical values of the test statistics need to be computed only once. The critical values of $\kappa_n$ and $\omega_n^2$ at commonly used significance levels are given in Table \ref{crit} below.

\begin{table}[t]
\vspace{8pt}
\begin{tabular}{c||c|c|c|c|}
  & Clayton & Clayton(2) & Gumbel & Gumbel(2)\\
  \hline \hline
  $\kappa_n$ & 41 & 45 & 46 & 46 \\
  \hline
  $\omega_n^2$ & 53 & 53 & 43 & 47 \\
  \hline
\end{tabular}\\
\vspace{10pt}
\begin{tabular}{c||c|c|c|c|}
  & Clayton & Clayton(2) & Gumbel & Gumbel(2)\\
  \hline \hline
  $\kappa_n$ & 11 & 11 & 6 & 14 \\
  \hline
  $\omega_n^2$ & 11 & 14 & 5 & 8 \\
  \hline
\end{tabular}
\vspace{10pt}
\caption{Number of rejections for 1000 samples at 5\% (top) and 1\% (bottom) significance levels under the various null hypotheses}
\label{size}
\end{table}
\begin{table}[t]
\begin{tabular}{c||c|c|c|}
  & 10\% & 5\% & 1\%\\
  \hline \hline
  $\kappa_n$ & 2.100 & 2.362 & 2.865 \\
  \hline
  $\omega_n^2$ & 0.526 & 0.708 & 1.186 \\
  \hline
\end{tabular}
\vspace{10pt}
\caption{Critical values of $\kappa_n$ and $\omega_n^2$ at various significance levels}
\vspace{-10pt}
\label{crit}
\end{table}

To observe the behavior of the test statistics under the alternative hypothesis, we generate 1000 out-of-model samples (of size $n=200$ as before) for each of the four models considered above. For the Clayton and Clayton(2) models, we generate samples from the bivariate distribution with Exponential(1) margins and Gumbel(2) copula. Note that the Gumbel(2) copula has the same Kendall's tau as Clayton(2), namely 1/2. For the Gumbel and Gumbel(2) models, we generate samples from the bivariate distribution with Lomax(3,1) margins and Clayton(2) copula. From each sample, we construct the test process $W_n$ on the grid $\mathcal{G}$ and compute the two test statistics $\kappa_n$ and $\omega_n^2$ as before. The resulting rejection frequencies at 5\% significance level are shown in Table \ref{power} below. These numbers confirm that tests based on our approach have high power, even with a moderate sample size of 200.

\begin{table}[h]
\label{power}
\begin{tabular}{r||c|c|c|c|}
  testing for: & Clayton & Clayton(2) & Gumbel & Gumbel(2) \\
  sampled from: & Gumbel(2) & Gumbel(2) & Clayton(2) & Clayton(2) \\
  \hline \hline
  $\kappa_n$ & 787 & 1000 & 956 & 768 \\
  \hline
  $\omega_n^2$ & 857 & 999 & 938 & 761 \\
  \hline
\end{tabular}
\vspace{10pt}
\caption{Number of rejections for 1000 samples at the 5\% significance level under the various alternative hypotheses}
\end{table}
\vspace{-15pt}

\subsection{Data analysis} We consider a data set consisting of log-concentra\-tions of seven metallic elements (uranium [U], lithium [Li], cobalt [Co], potassium [K], caesium [Cs], scandium [Sc], titanium [Ti]) in 655 water samples collected near Grand Junction, Colorado in the late 1970s. In \cite{cook:johnson:1986}, the pairwise dependence structures of U-Cs, Co-Ti and Cs-Sc log-concentrations were investigated, and it was found that the Clayton copula, or a two-parameter extension of it, provides a better fit (in terms of likelihood values) to each of these pairs than the normal copula, under the assumption of normal marginal distributions. This two-parameter extension of the Clayton copula is defined as:
\begin{align}\label{cj.cop}
\lefteqn{C_{\lambda_1, \lambda_2}(u,v)} \qquad \notag\\
&= (1 + \lambda_2)(u^{-\lambda_1} + v^{-\lambda_1} - 1)^{-1/\lambda_1} + \lambda_2(2u^{-\lambda_1} + 2v^{-\lambda_1} - 3)^{-1/\lambda_1}\\
& \hspace{26pt} {} - \lambda_2(2u^{-\lambda_1} + v^{-\lambda_1} - 2)^{-1/\lambda_1} - \lambda_2(u^{-\lambda_1} + 2v^{-\lambda_1} - 2)^{-1/\lambda_1}, \notag
\end{align}
for parameters $\lambda_1 > 0$ and $\lambda_2 \in [0,1]$. Note that when $\lambda_2=0$, the expression (\ref{cj.cop}) reduces to the usual Clayton copula with parameter $\lambda_1$.

For our analysis, we focus on the pair Co-Sc (which was not investigated in \cite{cook:johnson:1986}), since the assumption of normal margins seems most plausible for the Co and Sc log-concentrations; see normal QQ-plots in Fig.\ \ref{norm.qq} below. Also see Fig.\ \ref{scatter} for a scatter plot of the Co-Sc log-concentrations as well as a scatter plot of the rank-transformed data.

\begin{figure}[h]
\epsfig{file=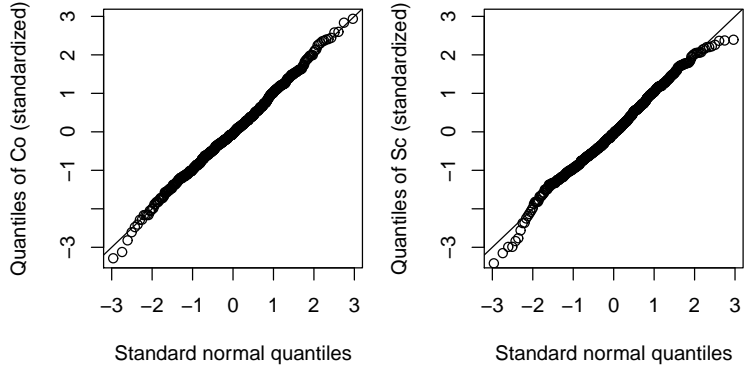,width=0.8\textwidth}
\caption{Normal QQ-plots for Co and Sc log-concentrations}
\label{norm.qq}
\end{figure}

\begin{figure}[h]
\epsfig{file=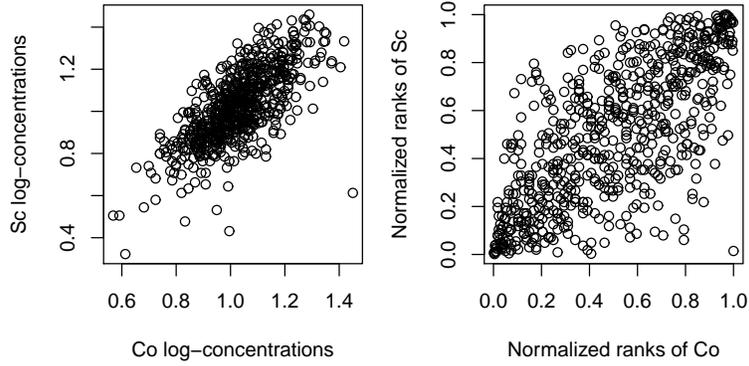,width=0.8\textwidth}
\caption{Scatter plots for Co and Sc log-concentrations (left) and normalized ranks of Co and Sc log-concentrations (right)}
\label{scatter}
\end{figure}

Under the assumption of normal margins, we test for four parametric copula families: Clayton, Frank, Gumbel, and the two-parameter family described in (\ref{cj.cop}), which we will call the \emph{Cook-Johnson copula}. Recall that the bivariate Frank family of copulas is given by
$$
C_\lambda(u,v) = -\frac{1}{\lambda} \log \bigg( 1 + \frac{(e^{-\lambda u} - 1)(e^{-\lambda v} - 1)}{e^{-\lambda} - 1} \bigg),
$$
for $\lambda \in \mathbb{R}$, with $C_0(u,v) = uv$. Following the methodology described in Subsection \ref{sim.stu}, we compute the test process $W_n$ and the corresponding test statistics $\kappa_n$ and $\omega_n^2$ for each of these four models. For the observed values of the test statistics, we compute $p$-values from the benchmark distribution tables constructed from the true bivariate standard Wiener process. The results are given in Table \ref{p.val}. We observe that all models except Frank are rejected at 5\% significance level, and Clayton and Cook-Johnson models in particular are rejected very strongly.

\begin{table}[h]
\vspace{8pt}
\begin{tabular}{c||c|c|c|c|c|}
  & Clayton & Frank & Gumbel & C-J\\
  \hline \hline
  $\kappa_n$ & 0.0000 & 0.1664 & 0.0318 & 0.0000\\
  \hline
  $\omega_n^2$ & 0.0000 & 0.1281 & 0.0278 & 0.0000\\
  \hline
\end{tabular}
\vspace{10pt}
\caption{$p$-values for various copula models for the Co-Sc log-concentrations, under assumption of marginal normality}
\vspace{-10pt}
\label{p.val}
\end{table}

\begin{figure}[t]
\epsfig{file=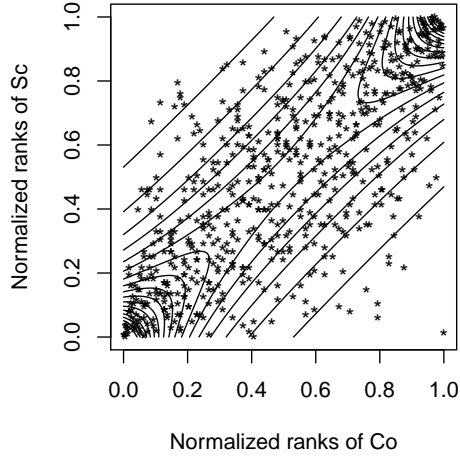,width=0.5\textwidth}
\caption{Contour plot of the fitted Frank copula density, superimposed on the scatter plot of the rank-transformed Co and Sc log-concentrations}
\label{contour}
\end{figure}

The model with normal margins and Frank copula yields maximum likelihood estimates of $\widehat{\mu}_{\text{Co}} = 1.025$ and $\widehat{\sigma}_{\text{Co}} = 0.136$ for the mean and standard deviation of Co log-concentrations, $\widehat{\mu}_{\text{Sc}} = 1.021$ and $\widehat{\sigma}_{\text{Sc}} = 0.178$ for the mean and standard deviation of Sc log-concentrations, and $\widehat{\lambda} = 6.589$ for the Frank copula parameter. The estimated value of $\lambda$ suggests moderate positive dependence, corresponding to a Kendall's $\tau$ of 0.544 and Spearman's $\rho$ of 0.743. Direct sample estimates for these coefficients are $0.535$ and $0.718$, respectively. The contour plot of the fitted Frank copula density, together with the scatter plot of the rank-transformed data, can be seen in Fig.\ \ref{contour}.

\section{Proofs}\label{proofs}

The proofs of Theorems \ref{thm.1} and \ref{thm.2} can be found in a supplementary document. We present the proofs of Theorems \ref{thm.3} and \ref{thm.4} below.

\subsection*{Proof of Theorem \ref{thm.3}} By Theorem \ref{thm.1} and Skorohod's representation theorem (\cite{billingsley:1999}, Theorem 6.7), there is a probability space where probabilistically equivalent versions of $\eta_n$ and $\eta$ are defined, and these satisfy $\| \eta_n - \eta \|_{[\delta,1-\delta]^d} \to 0$ a.s., with $\| \cdot \|_S := \sup_S| \cdot |$ for $S \subset [0,1]^d$. We will show that in this probability space,
\begin{equation}\label{show.1}
\| W_n - W \|_{[0,1]^d} \stackrel{P}{\to} 0,
\end{equation}
with $W$ as defined in (\ref{W.def}). In view of Theorem \ref{thm.2}, this will suffice for the proof.

Throughout the proof, we will let $A_\delta(\mathbf{u})$ denote the set $[\boldsymbol{\delta}, \boldsymbol{\delta}+(1-2\delta)\mathbf{u}]$ for $\mathbf{u} \in [0,1]^d$. Note that (\ref{show.1}) will follow from
\begin{equation}\label{show.2}
\bigg\| \int_{A_\delta(\mathbf{u})} \frac{1}{\sqrt{c_{\widehat{\boldsymbol{\lambda}}}(\mathbf{s})}} \, \dd \eta_n(\mathbf{s}) - \int_{A_\delta(\mathbf{u})} \frac{1}{\sqrt{c(\mathbf{s})}} \, \dd \eta(\mathbf{s}) \bigg\|_{[0,1]^d} \stackrel{P}{\to} 0
\end{equation}
and
\begin{equation}\label{show.3}
\begin{split}
\lefteqn{\bigg\| \int_{A_\delta(\mathbf{u})} \widehat{\mathbf{k}}(\mathbf{s})^\mathsf{T}
\bigg( \widehat{\mathbf{I}}_\delta^{-1}(s_d) \int_{S_\delta(s_d)} \widehat{\mathbf{k}}(\mathbf{s'})\,\dd \eta_n(\mathbf{s'}) \bigg) \sqrt{c_{\widehat{\boldsymbol{\lambda}}}(\mathbf{s})} \, \dd \mathbf{s}} \qquad \\
& - \int_{A_\delta(\mathbf{u})} \mathbf{k}(\mathbf{s})^\mathsf{T}
\bigg( \mathbf{I}_\delta^{-1}(s_d) \int_{S_\delta(s_d)} \mathbf{k}(\mathbf{s'})\,\dd \eta(\mathbf{s'}) \bigg) \sqrt{c(\mathbf{s})} \, \dd \mathbf{s}  \bigg\|_{[0,1]^d} \stackrel{P}{\to} 0.
\end{split}
\end{equation}
We will prove (\ref{show.2}) first. Let $\Delta_n := \eta_n - \eta$. Then (\ref{show.2}) will follow from
\begin{equation}\label{show.4}
\bigg\| \int_{A_\delta(\mathbf{u})} \Delta\gamma(\mathbf{s}) \,\dd \eta(\mathbf{s}) \bigg\|_{[0,1]^d} \stackrel{P}{\to} 0, \qquad \bigg\| \int_{A_\delta(\mathbf{u})} \widehat{\gamma}(\mathbf{s}) \,\dd \Delta_n(\mathbf{s}) \bigg\|_{[0,1]^d} \stackrel{P}{\to} 0.
\end{equation}
Applying integration by parts (\cite{henstock:1973}, Theorem 3) to the first integral term in (\ref{show.4}), we obtain the following bound:
\begin{equation}\label{ibp.1}
\begin{split}
\bigg| \int_{A_\delta(\mathbf{u})} \Delta\gamma(\mathbf{s}) \,\dd \eta(\mathbf{s}) \bigg|
&\le \sum_{\mathbf{v} \in \mathcal{V}_{A_\delta(\mathbf{u})}} \big| \Delta\gamma(\mathbf{v})\eta(\mathbf{v}) \big| + \| \eta \|_{A_\delta(\mathbf{u})} V_{A_\delta(\mathbf{u})}^{\text{HK}} (\Delta \gamma)\\
&\le \| \eta \|_\delta \big(2^d \| \Delta \gamma \|_\delta + V_\delta^{\text{HK}} (\Delta \gamma)\big),
\end{split}
\end{equation}
where $\mathcal{V}_{A_\delta(\mathbf{u})}$ denotes the set of the $2^d$ vertices of the hyperrectangle $A_\delta(\mathbf{u})$, and $\| \cdot \|_\delta$ is short-hand notation for $\| \cdot \|_{[\delta,1-\delta]^d}$. Now, Assumptions A2-A3 ensure that $\eta$ is continuous (hence bounded) on $[\delta,1-\delta]^d$, A4 ensures that $|\Delta \gamma|$ is $o_P(1)$ uniformly over $[\delta,1-\delta]^d$, and A6 ensures that $V_\delta^{\text{HK}} (\Delta \gamma)$ is $o_P(1)$ as well. It follows that the far right-hand side of (\ref{ibp.1}) vanishes in probability, and the first convergence in (\ref{show.4}) is proved. The second convergence in (\ref{show.4}) follows from a similar integration by parts argument:
\begin{equation*}
\bigg| \int_{A_\delta(\mathbf{u})} \widehat{\gamma}(\mathbf{s}) \,\dd \Delta_n(\mathbf{s}) \bigg| \le \| \Delta_n \|_\delta \big(2^d \| \widehat{\gamma} \|_\delta + V_\delta^{\text{HK}} (\widehat{\gamma})\big),
\end{equation*}
where the right-hand side is $o_P(1)$ since $\| \Delta_n \|_\delta$ is $o_P(1)$ and $\| \widehat{\gamma} \|_\delta$ as well as $V_\delta^{\text{HK}} (\widehat{\gamma})$ are $O_P(1)$ terms.

We have thus established (\ref{show.2}), and it remains to prove (\ref{show.3}). For ease of notation, we let
\begin{align*}
H(\mathbf{s}) &= \mathbf{k}(\mathbf{s})^\mathsf{T} \mathbf{I}_\delta^{-1}(s_d) \int_{S_\delta(s_d)} \mathbf{k}(\mathbf{s'})\,\dd \eta(\mathbf{s'}),\\
H_n(\mathbf{s}) &= \mathbf{k}(\mathbf{s})^\mathsf{T} \mathbf{I}_\delta^{-1}(s_d) \int_{S_\delta(s_d)} \mathbf{k}(\mathbf{s'})\,\dd \eta_n(\mathbf{s'}),\\
\widehat{H}(\mathbf{s}) &= \widehat{\mathbf{k}}(\mathbf{s})^\mathsf{T} \widehat{\mathbf{I}}_\delta^{-1}(s_d) \int_{S_\delta(s_d)} \widehat{\mathbf{k}}(\mathbf{s'})\,\dd \eta(\mathbf{s'}),\\
\widehat{H}_n(\mathbf{s}) &= \widehat{\mathbf{k}}(\mathbf{s})^\mathsf{T} \widehat{\mathbf{I}}_\delta^{-1}(s_d) \int_{S_\delta(s_d)} \widehat{\mathbf{k}}(\mathbf{s'})\,\dd \eta_n(\mathbf{s'}).
\end{align*}
Then (\ref{show.3}) can be written succinctly as
\begin{equation*}
\bigg\| \int_{A_\delta(\mathbf{u})} \Big(\widehat{H}_n(\mathbf{s}) \sqrt{c_{\widehat{\boldsymbol{\lambda}}}(\mathbf{s})} -  H(\mathbf{s}) \sqrt{c(\mathbf{s})} \Big) \, \dd \mathbf{s}  \bigg\|_{[0,1]^d} \stackrel{P}{\to} 0,
\end{equation*}
which can be proved by showing
\begin{equation}\label{show.5}
\big\| H \big(\sqrt{c_{\widehat{\boldsymbol{\lambda}}}} - \sqrt{c} \big) \big\|_\delta \stackrel{P}{\to} 0, \quad \big\| \big(\widehat{H}_n - H \big) \sqrt{c_{\widehat{\boldsymbol{\lambda}}}} \big\|_\delta \stackrel{P}{\to} 0.
\end{equation}
The first convergence in (\ref{show.5}) follows easily from the continuity (hence boundedness) of $H$ over $[\delta,1-\delta]^d$ and the continuity of $\sqrt{c_{\boldsymbol{\lambda}}(\mathbf{u})}$ over $(\mathbf{u},\boldsymbol{\lambda}) \in [\delta, 1-\delta]^d \times \Lambda$. As for the second convergence in (\ref{show.5}), since $\big\| \sqrt{c_{\widehat{\boldsymbol{\lambda}}}} \big\|_\delta = O_P(1)$, we need to show that $\| \widehat{H}_n - H \|_\delta \stackrel{P}{\to} 0$. We will do this by proving
\begin{equation}\label{show.6}
\| H_n - H \|_\delta \stackrel{P}{\to} 0, \quad \| \widehat{H}_n - H_n \|_\delta \stackrel{P}{\to} 0.
\end{equation}
Consider the first convergence in (\ref{show.6}). We have
\begin{equation*}
\| H_n - H \|_\delta = \bigg\| \mathbf{k}(\mathbf{s})^\mathsf{T} \mathbf{I}_\delta^{-1}(s_d) \int_{S_\delta(s_d)}  \mathbf{k}(\mathbf{s'})\,\dd \Delta_n(\mathbf{s'}) \bigg\|_\delta,
\end{equation*}
with $\Delta_n = \eta_n-\eta$, as before. The term $|\mathbf{k}(\mathbf{s})^\mathsf{T} \mathbf{I}_\delta^{-1}(s_d)|$ is component-wise bounded on $[\delta,1-\delta]^d$ by continuity, so we need to show that
\begin{equation}\label{show.7}
\sup_{t \in [\delta,1-\delta]} \bigg| \int_{S_\delta(t)} k_i(\mathbf{s'})\,\dd \Delta_n(\mathbf{s'}) \bigg| \stackrel{P}{\to} 0, \quad i=1,\ldots, 1+dm+p.
\end{equation}
Applying integration by parts as before, we obtain
\begin{equation*}
\bigg| \int_{S_\delta(t)} k_i(\mathbf{s'})\,\dd \Delta_n(\mathbf{s'}) \bigg| \le \| \Delta_n \|_\delta \big(2^d \| k_i \|_\delta + V_\delta^{\text{HK}} (k_i)\big),
\end{equation*}
where the right-hand side is $o_P(1)$ since $\| k_i \|_\delta < \infty$, $V_\delta^{\text{HK}} (k_i) < \infty$ and $\| \Delta_n \|_\delta = o_P(1)$. Hence (\ref{show.7}) is established and it remains to prove the second convergence in (\ref{show.6}).

By virtue of the first convergence in (\ref{show.6}), and an analogous result for $\widehat{H}_n$ and $\widehat{H}$, it will suffice to prove $\| \widehat{H} - H \|_\delta \stackrel{P}{\to} 0$. Note that
\begin{equation}\label{note.1}
\begin{split}
|\widehat{H}(\mathbf{s}) - H(\mathbf{s})|
&\le \big| \widehat{\mathbf{k}}(\mathbf{s})^\mathsf{T} \widehat{\mathbf{I}}_\delta^{-1}(s_d) - \mathbf{k}(\mathbf{s})^\mathsf{T} \mathbf{I}_\delta^{-1}(s_d) \big| \cdot \bigg| \int_{S_\delta(s_d)} \mathbf{k}(\mathbf{s'})\,\dd \eta(\mathbf{s'}) \bigg|\\
& \qquad + \big| \widehat{\mathbf{k}}(\mathbf{s})^\mathsf{T} \widehat{\mathbf{I}}_\delta^{-1}(s_d) \big| \cdot \bigg| \int_{S_\delta(s_d)} \big(\widehat{\mathbf{k}}(\mathbf{s}') - \mathbf{k}(\mathbf{s'})\big) \,\dd \eta(\mathbf{s'}) \bigg|,
\end{split}
\end{equation}
where absolute values should be interpreted component-wise, as usual. Consider the first term on the right-hand side of (\ref{note.1}). Since the mapping
\begin{equation*}
(\mathbf{s}, \boldsymbol{\theta}'_1, \ldots, \boldsymbol{\theta}'_d, \boldsymbol{\lambda}') \mapsto \mathbf{k}(\mathbf{s}, \boldsymbol{\theta}'_1, \ldots, \boldsymbol{\theta}'_d, \boldsymbol{\lambda}')
\end{equation*}
is continuous over $[\delta,1-\delta]^d \times \Theta^d \times \Lambda$, the difference
\begin{equation*}
\big| \widehat{\mathbf{k}}(\mathbf{s})^\mathsf{T} \widehat{\mathbf{I}}_\delta^{-1}(s_d) - \mathbf{k}(\mathbf{s})^\mathsf{T} \mathbf{I}_\delta^{-1}(s_d) \big|
\end{equation*}
is $o_P(1)$ uniformly over $\mathbf{s} \in [\delta,1-\delta]^d$. Moreover, an integration by parts argument as before yields that
\begin{equation*}
\bigg| \int_{S_\delta(s_d)} k_i(\mathbf{s'})\,\dd \eta(\mathbf{s'}) \bigg| \le \| \eta \|_\delta \big(2^d\| k_i \|_\delta + V^{\text{HK}}_\delta(k_i)\big),
\end{equation*}
for $i=1,\ldots, 1+dm+p$, where the right-hand side is $O_P(1)$. So the first summand on the right-hand side of (\ref{note.1}) is $o_P(1)$ uniformly over $\mathbf{s} \in [\delta,1-\delta]^d$. The second summand there can be handled similarly: the term $\big| \widehat{\mathbf{k}}(\mathbf{s})^\mathsf{T} \widehat{\mathbf{I}}_\delta^{-1}(s_d) \big|$ is $O_P(1)$, and integration by parts yields
\begin{equation*}
\bigg| \int_{S_\delta(s_d)} \Delta k_i(\mathbf{s'})\,\dd \eta(\mathbf{s'}) \bigg| \le \| \eta \|_\delta \big(2^d \| \Delta k_i \|_\delta + V^{\text{HK}}_\delta(\Delta k_i)\big)
\end{equation*}
for $i=1,\ldots,1+dm+p$, where the right-hand side is $o_P(1)$.

Both convergences in (\ref{show.5}) are thereby established, which in turn proves (\ref{show.3}). \hfill $\Box$

\subsection*{Proof of Theorem \ref{thm.4}} We have, from the Cameron-Martin-Girsanov theorem, that
\begin{align*}
\log\Big(\frac{\dd \widetilde{Q}}{\dd Q}\Big) &= (1-2\delta)^{d/2} \int_{[0,1]^d} g(\boldsymbol{\delta}+(1-2\delta)\mathbf{u}) \sqrt{c_{\boldsymbol{\lambda}_0}(\boldsymbol{\delta}+(1-2\delta)\mathbf{u})} \, \dd V(\mathbf{u})\\
& \qquad {} - \frac{(1-2\delta)^d}{2} \int_{[0,1]^d} g^2(\boldsymbol{\delta}+(1-2\delta)\mathbf{u}) \, c_{\boldsymbol{\lambda}_0}(\boldsymbol{\delta}+(1-2\delta)\mathbf{u}) \, \dd \mathbf{u},
\end{align*}
with $V = W$ under $Q$ and $V = \widetilde{W}$ under $\widetilde{Q}$, and $g$ as defined in (\ref{g.def}). This immediately yields
\begin{equation*}
\log\Big(\frac{\dd\widetilde{Q}}{\dd Q}\Big) \sim 
\begin{cases} 
N\big(-\frac{1}{2} \|g\|^2, \|g\|^2\big) &\text{ under } Q,\\[5pt]
N\big(\frac{1}{2} \|g\|^2, \|g\|^2\big) &\text{ under } \widetilde{Q},
\end{cases}
\end{equation*}
with $\|g\|$ as in (\ref{add.3}), with the convention that $g=0$ outside of $[\delta,1-\delta]^d$. Using the fact that
\begin{equation*}
d(\widetilde{Q},Q) = \widetilde{Q}\big[\log(\dd\widetilde{Q}/\dd Q) > 0\big] - Q\big[\log(\dd\widetilde{Q}/\dd Q) > 0\big],
\end{equation*} 
we also obtain $d(\widetilde{Q},Q) = \nu(g) = 2\Phi(\frac{1}{2} \|g\|) - 1$.

Thus it remains to show that $\|g\| = \|h\|$, which is equivalent to
\begin{equation}\label{toshow2}
\begin{split}
\lefteqn{2\int_{[\delta,1-\delta]^d} \mathbf{k}(\mathbf{s})^\mathsf{T} \bigg( \mathbf{I}_\delta^{-1}(s_d) \int_{S_\delta(s_d)} \mathbf{k}(\mathbf{s'})h(\mathbf{s'})\,\dd C_{\boldsymbol{\lambda}_0}(\mathbf{s'}) \bigg) h(\mathbf{s}) \, \dd C_{\boldsymbol{\lambda}_0}(\mathbf{s})} \hspace{10pt}\\
&= \int_{[\delta,1-\delta]^d} \bigg[ \mathbf{k}(\mathbf{s})^\mathsf{T} \bigg( \mathbf{I}_\delta^{-1}(s_d) \int_{S_\delta(s_d)} \mathbf{k}(\mathbf{s'})h(\mathbf{s'})\,\dd C_{\boldsymbol{\lambda}_0}(\mathbf{s'}) \bigg) \bigg]^2 \, \dd C_{\boldsymbol{\lambda}_0}(\mathbf{s}).
\end{split}
\end{equation}
For ease of notation, let $E_1$ and $E_2$ denote the left- and right-hand sides of (\ref{toshow2}), respectively. Also let $S^-_\delta(t) = [\delta,1-\delta/2]^d \setminus S_\delta(t)$ for $t \in [\delta,1-\delta/2)$ and define
\begin{equation*}
\mathbf{H}(t) = \int_{S^-_\delta(t)} \mathbf{k}(\mathbf{s}) h(\mathbf{s})\, \dd C_{\boldsymbol{\lambda}_0}(\mathbf{s}), \quad t \in [\delta,1-\delta/2).
\end{equation*}
We have
\begin{align*}
E_1
&= 2\int_{[\delta,1-\delta]^d} \bigg[ \int_{S_\delta(s_d)} \mathbf{k}(\mathbf{s'})^\mathsf{T}  h(\mathbf{s'}) \,\dd C_{\boldsymbol{\lambda}_0}(\mathbf{s'}) \bigg] \mathbf{I}_\delta^{-1}(s_d) \mathbf{k}(\mathbf{s})h(\mathbf{s}) \, \dd C_{\boldsymbol{\lambda}_0}(\mathbf{s})\\
&= -2\int_{[\delta,1-\delta]^d} \bigg[ \int_{S^-_\delta(s_d)} \mathbf{k}(\mathbf{s'})^\mathsf{T}  h(\mathbf{s'}) \,\dd C_{\boldsymbol{\lambda}_0}(\mathbf{s'}) \bigg] \mathbf{I}_\delta^{-1}(s_d) \mathbf{k}(\mathbf{s})h(\mathbf{s}) \, \dd C_{\boldsymbol{\lambda}_0}(\mathbf{s})\\
&= -2\int_\delta^{1-\delta} \bigg[ \int_{S^-_\delta(s_d)} \mathbf{k}(\mathbf{s'})^\mathsf{T}  h(\mathbf{s'}) \,\dd C_{\boldsymbol{\lambda}_0}(\mathbf{s'}) \bigg] \mathbf{I}_\delta^{-1}(s_d)\\
& \hspace{100pt} {} \times \bigg[ \int_{[\delta,1-\delta]^{d-1}}\mathbf{k}(\mathbf{s})h(\mathbf{s})c_{\boldsymbol{\lambda}_0}(\mathbf{s}) \, \dd s_1 \ldots \dd s_{d-1} \bigg] \dd s_d\\
&= -2 \int_\delta^{1-\delta} \mathbf{H}(s_d)^\mathsf{T} \mathbf{I}_\delta^{-1}(s_d) \, \dd \mathbf{H}(s_d),
\end{align*}
where the second equality above follows from
\begin{equation*}
\int_{[\delta,1-\delta/2]^d} \mathbf{k}(\mathbf{s}) h(\mathbf{s})\, \dd C_{\boldsymbol{\lambda}_0}(\mathbf{s}) = \mathbf{0},
\end{equation*}
which is a consequence of Assumption B0. Now, denoting $\mathbf{G}(t) = \mathbf{I}_\delta^{-1}(t)\mathbf{H}(t)$ and applying integration by parts, we obtain
\begin{equation*}
E_1 = -2 \int_\delta^{1-\delta} \mathbf{G}(s_d)^\mathsf{T} \, \dd \mathbf{H}(s_d) = 2 \int_\delta^{1-\delta} \mathbf{H}(s_d)^\mathsf{T} \, \dd \mathbf{G}(s_d).
\end{equation*}
By the product rule of differentiation, we can write
\begin{equation*}
\dd \mathbf{G}(s_d) = \big[(\mathbf{I}_\delta^{-1})'(s_d)\mathbf{H}(s_d) + \mathbf{I}_\delta^{-1}(s_d)\mathbf{H}'(s_d)\big] \, \dd s_d,
\end{equation*}
where derivatives should be interpreted component-wise. Using the identity
\begin{equation*}
(\mathbf{I}_\delta^{-1})'(s_d) = - \mathbf{I}_\delta^{-1}(s_d) \mathbf{I}'_\delta(s_d)\mathbf{I}_\delta^{-1}(s_d),
\end{equation*}
we obtain
\begin{align*}
E_1 &= -2 \int_\delta^{1-\delta} \mathbf{H}(s_d)^\mathsf{T} \mathbf{I}_\delta^{-1}(s_d) \mathbf{I}'_\delta(s_d)\mathbf{I}_\delta^{-1}(s_d) \mathbf{H}(s_d) \, \dd s_d\\
& \hspace{130pt} {} + 2 \int_\delta^{1-\delta} \mathbf{H}(s_d)^\mathsf{T} \mathbf{I}_\delta^{-1}(s_d) \mathbf{H}'(s_d) \, \dd s_d\\
&= -2 \int_\delta^{1-\delta} \mathbf{H}(s_d)^\mathsf{T} \mathbf{I}_\delta^{-1}(s_d) \mathbf{I}'_\delta(s_d)\mathbf{I}_\delta^{-1}(s_d) \mathbf{H}(s_d) \, \dd s_d - E_1.
\end{align*}
It follows that
\begin{align*}
E_1 &= - \int_\delta^{1-\delta} \mathbf{H}(s_d)^\mathsf{T} \mathbf{I}_\delta^{-1}(s_d) \mathbf{I}'_\delta(s_d)\mathbf{I}_\delta^{-1}(s_d) \mathbf{H}(s_d) \, \dd s_d\\
&= \int_\delta^{1-\delta} \mathbf{H}(s_d)^\mathsf{T} \mathbf{I}_\delta^{-1}(s_d) \bigg[\int_{[\delta,1-\delta/2]^{d-1}} \mathbf{k}(\mathbf{s})\mathbf{k}(\mathbf{s})^\mathsf{T}c_{\boldsymbol{\lambda}_0}(\mathbf{s})\dd s_1 \ldots \dd s_{d-1}\bigg]\\ & \hspace{230pt} {} \cdot \mathbf{I}_\delta^{-1}(s_d) \mathbf{H}(s_d) \, \dd s_d\\
&= \int_{[\delta,1-\delta]^d} \mathbf{H}(s_d)^\mathsf{T} \mathbf{I}_\delta^{-1}(s_d) \mathbf{k}(\mathbf{s})\mathbf{k}(\mathbf{s})^\mathsf{T}\mathbf{I}_\delta^{-1}(s_d) \mathbf{H}(s_d) \dd C_{\boldsymbol{\lambda}_0}(\mathbf{s})\\
&= \int_{[\delta,1-\delta]^d} \Big[ \mathbf{k}(\mathbf{s})^\mathsf{T}\mathbf{I}_\delta^{-1}(s_d) \mathbf{H}(s_d)\Big]^2 \dd C_{\boldsymbol{\lambda}_0}(\mathbf{s})\\
&= E_2,
\end{align*}
as desired. Thus (\ref{toshow2}) is established. \hfill $\Box$

\section*{Acknowledgements}

We are very grateful to Estate V.\ Khmaladze for stimulating discussions about this paper. We are also grateful to the participants of the Eurandom-ISI Workshop on Actuarial and Financial Statistics, the $10^\text{th}$ Extreme Value Analysis Conference, the $4^\text{th}$ Conference of the International Society for Nonparametric Statistics, as well as seminar participants at Tilburg and Bocconi Universities, for their comments and suggestions.

\bibliographystyle{apalike}
\bibliography{copula}
\end{document}